\magnification=1000
\hsize=11.7cm
\vsize=18.9cm
\lineskip2pt \lineskiplimit2pt
\nopagenumbers

\hoffset=-1truein
\voffset=-1truein

\advance\voffset by 4truecm
\advance\hoffset by 4.5truecm

\newif\ifentete

\headline{\ifentete\ifodd	\count0 
      \rlap{\head}\hfill\tenrm\llap{\the\count0}\relax
    \else
        \tenrm\rlap{\the\count0}\hfill\llap{\head} \relax
    \fi\else
\global\entetetrue\fi}

\def\entete#1{\entetefalse\gdef\head{#1}}
\entete{}

\input amssym.def
\input amssym.tex

\def\-{\hbox{-}}
\def\.{{\cdot}}
\def\O{{\cal O}}
\def\K{{\cal K}}
\def\F{{\cal F}}

\def\G{{\cal G}}

\def\ch{\frak c\frak h}

\def\Gr{\frak G\frak r}

\def\int{\frak i\frak n\frak t}

\def\qq{\quad{\rm and}\quad}

\def\aut{\frak a\frak u\frak t}

 3
 2
\font\large=cmr10  scaled \magstep 2
 2
 2
 2
\font\cds=cmr7

\centerline{\large Note on the reduction of Alperin's Conjecture}

\vskip 0.5cm

\centerline{\bf Lluis Puig }
\medskip
\centerline{\cds CNRS, Institut de Math\'ematiques de Jussieu}
\smallskip
\centerline{\cds 6 Av Bizet, 94340 Joinville-le-Pont, France}
\smallskip
\centerline{\cds puig@math.jussieu.fr}

\vskip 0.5cm
£1. In a recent paper [2], Gabriel Navarro and Pham Huu Tiep show that the so-called Alperin Weight Conjecture
can be verified {\it via\/}  the {\it Classification of the Finite Simple Groups\/}, provided any simple group fulfills
a very precise list of conditions that they consider easier to check than ours, firstly stated in [3,~Theorem~16.45] and
significantly weakened in [4,~Theorem~1.6]{\footnote{\dag}{\cds  Gabriel Navarro and Pham Huu Tiep pointed out to us that, when submitting [2], they were not aware of our paper [4], 
only available in arXiv since April 2010.}}.

\medskip
£2. Actually, in the introduction of [3] --- from I29 to I37 --- we consider the most precise  Alperin's Conjecture concerning
any block of a finite group, and introduce a refinement to this conjecture; but, only in~[4] we really show that its verification can be reduced to check that the {\it same\/} refinement holds on the so-called {\it quasi-simple\/} groups. To carry out this checking obviously depends on admitting the 
{\it Classification of the Finite Simple Groups\/}, and our proof of the reduction itself uses the {\it solvability\/} of the
{\it outer automorphism group\/} of a finite simple group, a known fact whose actual proof depends on this classification.

\medskip
£3. Our purpose here is to show to the interested reader{\footnote{\dag\dag}{\cds Indeed, a first submission of [4]
has been rejected since, according to the referee, ``the results (...) are too specialized for Inventiones''.}} that the results in [3] and the reduction arguments in [4] suggest a numerical statement --- implying Alperin's Conjecture --- which can be reduced again to check that the same  holds on the  quasi-simple groups and, this time, this statement on the
quasi-simple groups follows from the list of conditions demanded in~[2].

\medskip
£4. Let us be more explicit. Let $p$ be a prime number, $k$ an algebraically closed  field of characteristic $p\,,$
 $\O$ a complete discrete valuation ring of characteristic zero admitting $k$ as the {\it residue\/} field, and 
$\K$ the field of fractions of~$\O\,.$ Moreover, let $\hat G$ be a $k^*\-$group of finite $k^*\-$quotient~$G$ [3,~1.23], $b$ a block of $\hat G$  [3,~1.25] and $\G_k (\hat G,b)$
 the {\it scalar extension\/} from $\Bbb Z$ to $\O$ of the {\it Grothendieck group\/} of the category of finitely generated $k_*\hat Gb\-$modules [3,~14.3]. In [3,~Chap.~14], choosing a maximal Brauer $(b,\hat G)\-$pair $(P,e)\,,$ the existence of a suitable $k^*\-\Gr\-$valued functor $\widehat\aut_{(\F_{\!(b,\hat G)})^{^{\rm nc}}}$ over some full subcategory $(\F_{\!(b,\hat G)})^{^{\rm nc}}$ of the  {\it Frobenius $P\-$category\/}~$\F_{\!(b,\hat G)}$ [3,~3.2] allows us to consider an inverse limit of Grothendieck groups --- noted  $\G_k (\F_{\!(b,\hat G)},\widehat\aut_{(\F_{\!(b,\hat G)})^{^{\rm nc}}})$ and called the {\it Grothendieck group of\/} $\F_{\!(b,\hat G)}$ --- such that Alperin's Conjecture is actually equivalent to the existence of an $\O\-$module isomorphism [3,~I32 and Corollary~14.32]
 $$\G_k (\hat G,b)\cong \G_k (\F_{\!(b,\hat G)},\widehat\aut_{(\F_{\!(b,\hat G)})^{^{\rm nc}}})
 \eqno £4.1.$$

 \medskip
 £5. Denote by ${\rm Out}_{k^*}(\hat G)$ the group of {\it outer\/} $k^*\-$automorphisms of $\hat G$ and by 
${\rm Out}_{k^*}(\hat G)_b$ the stabilizer of $b$ in ${\rm Out}_{k^*}(\hat G)\,;$ it is clear that 
${\rm Out}_{k^*}(\hat G)_b$ acts on $\G_k (\hat G,b)\,,$ and in [3,~16.3 and~16.4] we show that this group still acts on~$\G_k (\F_{\!(b,\hat G)},\widehat\aut_{(\F_{\!(b,\hat G)})^{^{\rm nc}}})\,.$ Denoting by 
$${}^\K\G_k (\hat G,b)\qq {}^\K\G_k (\F_{\!(b,\hat G)},\widehat\aut_{(\F_{\!(b,\hat G)})^{^{\rm nc}}})
\eqno £5.1\phantom{.}$$
the respective {\it scalar extensions\/} of $\G_k (\hat G,b)$ and $\G_k (\F_{\!(b,\hat G)},
\widehat\aut_{(\F_{\!(b,\hat G)})^{^{\rm nc}}})$ from $\O$ to $\K\,,$ here we replace the statement (Q) in [4,~1.4] by the following statement
\smallskip
\noindent
(${}^\K$Q)\quad {\it For any $k^*\-$group $\hat G$ with finite $k^*\-$quotient $G$ and any block $b$ of $\hat G\,,$ 
there is an $\K {\rm Out}_{k^*}(\hat G)_b\-$module isomorphism 
$${}^\K \G_k (\hat G,b)\cong {}^\K \G_k (\F_{\!(b,\hat G)},\widehat\aut_{(\F_{\!(b,\hat G)})^{^{\rm nc}}})
\eqno £5.2\phantom{.}$$\/}

\medskip
£6. The big difference between the statements (Q) and (${}^\K$Q) is that the second one is equivalent to the equality
of the corresponding {\it $\K\-$characters\/} of~${\rm Out}_{k^*}(\hat G)_b\,;$ moreover, since these 
$\K {\rm Out}_{k^*}(\hat G)_b\-$modules  actually come from $\Bbb Q {\rm Out}_{k^*}(\hat G)_b\-$modules, isomorphism~£5.2 is finally equivalent
to the equalities
$${\rm rank}_\O\big(\G_k (\hat G,b)^C\big) = {\rm rank}_\O \big(\G_k (\F_{\!(b,\hat G)},
\widehat\aut_{(\F_{\!(b,\hat G)})^{^{\rm nc}}})^C\big)
\eqno £6.1\phantom{.}$$
for any cyclic subgroup $C$ of ${\rm Out}_{k^*}(\hat G)_b\,.$ On the other hand, without any change, our proof of [4, Theorem~1.6] still proves the following result.

\bigskip
\noindent
{\bf Theorem~£7.} {\it Assume that any block $(c,\hat H)$ having a normal sub-block~$(d,\hat S)$ of positive defect 
such that the $k^*\-$quotient $S$ of $\hat S$ is simple, $H/S$ is a cyclic $p'\-$group and $C_{H}(S) = \{1\}\,,$ fulfills the following two conditions
\smallskip
\noindent
{\rm £7.1}\quad ${\rm Out}(S)$ is solvable.
\smallskip
\noindent
{\rm £7.2}\quad For any cyclic subgroup $C$ of ${\rm Out}_{k^*}(\hat H)_c$ we have
$${\rm rank}_\O\big(\G_k (\hat H,c)^C\big) = {\rm rank}_\O \big(\G_k (\F_{\!(c,\hat H)},
\widehat\aut_{(\F_{\!(c,\hat H)})^{^{\rm nc}}})^C\big)\,.$$
Then, for any block $(b,\hat G)$ there is an $\K {\rm Out}_{k^*} (\hat G)_b\-$module
isomorphism
$${}^\K \G_k (\hat G,b)\cong {}^\K \G_k(\F_{\!(b,\hat G)},\widehat \aut_{(\F_{\!(b,\hat G)})^{^{\rm nc}}})
\eqno £7.3.$$\/}
\eject

\medskip
£8. Morover, with the same notation of the theorem, it is easily checked that the proof of [3, Corollary~14.32]
actually also proves that 
$$\eqalign{{\rm rank}_\O \big(\G_k (\F_{\!(c,\hat H)},
&\widehat\aut_{(\F_{\!(c,\hat H)})^{^{\rm nc}}})^C\big)\cr
&= \sum_{(\frak q,\Delta_n)} (-1)^{n}\,{\rm rank}_{\O} \Big(\G_k \big(\hat\F_{\!(c,\hat H)} (\frak q)\big)^{C_{\frak q}}\Big)\cr}
\eqno £8.1\phantom{.}$$
where $(\frak q,\Delta_n)$ runs over a set of representatives for the set of isomorphism
classes of regular $\ch^*(\F^{^{\rm sc}}_{\!(c,\hat H)})\-$objects [3, 45.2] and, for such a $(\frak q,\Delta_n)\,,$
$C_\frak q$ denotes the ``stabilizer'' of $(\frak q,\Delta_n)$ in $C$ (see [3,~15.33] for a similar notation); indeed, in that proof, the sequence of $C\-$fixed elements in the exact sequence~[3,~14.32.4] remains exact, since we are working over 
$\K\,,$ and then equality~£8.1 follows easily.

\medskip
£9. At this point, the old argument of Reinhard Kn\"orr and Geoffrey Robinson in [1], suitably adapted, shows that,
when proving statement (${}^\K$Q) arguing by induction on $\vert G\vert\,,$ we may assume that
$$\eqalign{\sum_{(\frak q,\Delta_n)} (-1)^{n}\,{\rm rank}_{\O} \Big(&\G_k \big(\hat\F_{\!(c,\hat H)} 
(\frak q)\big)^{C_{\frak q}}\Big)\cr
& = \sum_{(Q,f)} {\rm rank}_{\O} \Big(\G_k \big(\hat\F_{\!(c,\hat H)} (Q,f),\bar b_f\big)^{C_f}\Big)\cr}
\eqno £9.1\phantom{.}$$
where $(Q,f)$ runs over a set of representatives for the set of $H\-$conjugacy classes of selfcentralizing 
Brauer $(c,\hat H)\-$pairs [3,~1.16 and Corollary~7.3] and, for such a  selfcentralizing 
Brauer $(c,\hat H)\-$pair $(Q,f)\,,$ we denote by $\bar b_f$ the sum of blocks of defect zero of 
$\hat\F_{\!(c,\hat H)} (Q,f)/\F_Q (Q)\,,$ and by $C_f$ the ``stabilizer'' of~$(Q,f)$ in $C\,.$

\medskip
£10. That is to say, according to equalities~£8.1 and £9.1, when proving Theorem~£7 arguing by induction on 
$\vert G\vert$ we may replace condition~£7.2 by the alternative condition:
\smallskip
\noindent
£10.1\quad {\it For any cyclic subgroup $C$ of ${\rm Out}_{k^*}(\hat H)_c$ we have
$${\rm rank}_\O\big(\G_k (\hat H,c)^C\big) = \sum_{(Q,f)} {\rm rank}_{\O} \Big(\G_k \big(\hat\F_{\!(c,\hat H)} 
(Q,f),\bar b_f\big)^{C_f}\Big)$$
where $(Q,f)$ runs over a set of representatives for the set of $H\-$conjugacy classes of selfcentralizing 
Brauer $(c,\hat H)\-$pairs.\/}
\smallskip
\noindent
Actually, the corresponding form of statement (${}^\K$Q) is nothing but the so-called Equivariant form
of Alperin's Conjecture, somewhere stated by Geoffrey Robinson.

\medskip
£11. Finally, we claim that this condition follows from the conditions in~[2,~\S3] and the ``compatibility''
admitted in [2, Remark~3.1]. Indeed, first of all note that, following the terminology in~[2], we are only concerned 
by the {\it $p\-$radical\/} subgroups $Q$ of $\hat S$ such that the quotient $\bar N_{\hat S}(Q) = N_{\hat S}(Q)/Q$ 
admits a block of defect zero or, equivalently, a projective simple module~$M\,;$ in this case, since the restriction of $M$ to
$\bar C_{\hat S}(Q) = C_{\hat S}(Q)/Z(Q)$ remains projective and semisimple, it involves a block $\bar f$
of $\bar C_{\hat S}(Q)$ of defect zero and therefore, denoting by $f$ the corresponding block of $C_{\hat S} (Q)\,,$
the Brauer $\hat S\-$pair $(Q,f)$ is {\it selfcentralizing\/} [3,~1.16 and Corollary~7.3]. Then, recalling that there is a bijection
between selfcentralizing Brauer pairs and selfcentralizing local pointed groups [3,~7.4], we have in~[3,~Lemma~15.16]
a precise relationship between the sets of selfcentralizing Brauer $\hat H\-$ and $\hat S\-$pairs.

\medskip
£12.  Set $A = H/S$ and, for any irreducible Brauer character $\theta$ of $\hat S$ in the block $d\,,$
respectively denote by $A_\theta$ and $\hat H_\theta$ the stabilizers of $\theta$ in $A$ and~$\hat H\,;$ 
denoting by $c_\theta$ the block of $\hat H_\theta$ determined by $\theta$ and by 
$\G_k (\hat H_\theta,c_\theta\!\mid\! \theta)$ the corresponding direct summand of $\G_k (\hat H_\theta,c_\theta)\,,$ 
it is quite clear that we have
$${\rm rank}_\O\big(\G_k (\hat H,c)^C\big) = \sum_\theta {\rm rank}_\O
\big(\G_k (\hat H_\theta,c_\theta\!\mid\! \theta)^{C_\theta}\big)
\eqno £12.1\phantom{.}$$
where $\theta$ runs over a set of representatives for the set of $H\-$orbits of irreducible Brauer character of $\hat S$ 
and, for such a $\theta\,,$ $C_\theta$ denotes the ``stabilizer'' of~$\theta$ in~$C\,.$ Moreover, forgetting  the block
$c_\theta$ and denoting by $\G_k (\hat H_\theta\!\mid\! \theta)$~the corresponding direct summand of 
$\G_k (\hat H_\theta)\,,$ it follows from the so-called {\it Clifford theory\/} that we have a canonical isomorphism
$$\G_k (\hat H_\theta\!\mid\! \theta)\cong \G_k (\hat A^{^\theta}_\theta)
\eqno £12.2\phantom{.}$$
for a suitable central $k^*\-$extension $\hat A^{^\theta}_\theta$ of $A_\theta\,;$ note that, since $A_\theta$
 is cyclic, the  $k^*\-$extension $\hat A^{^\theta}_\theta$ is split and, since ${\rm Aut} (A_\theta)$ is abelian,
this automorphism group acts {\it canonically\/} on~$\hat A^{^\theta}_\theta\,.$

\medskip
£13. But, according to the conditions in~[2] and to our remark above, $\theta$ determines up to $S\-$conjugation
a selfcentralizing Brauer $(d,\hat S)\-$pair $(Q,f)$ together with a projective simple $k_*\bar N_{\hat S} (Q,f)\bar f\-$module
$M^*$ which necessarily has the form
$$M^*\cong W\otimes_k \bar M^*
\eqno £13.1\phantom{.}$$
where $W$ is a projective simple $k_* \bar C_{\hat S}(Q)\bar f\-$module, suitable extended to the corresponding central
$k^*\-$extension of $\bar N_S (Q,f)\,,$ and $\bar M^*$ is a projective simple $k_*\hat \F_{\!(d,\hat S)} (Q,f)\-$module,
considered as a module over  the corresponding central
$k^*\-$extension of $\bar N_S (Q,f)\,;$ let us respectively denote by $\theta^*$ and $\bar\theta^*$ the (Brauer) characters of $M^*$ and~$\bar M^*\,.$
\eject

\medskip
£14. Always according to the conditions in~[2], the uniqueness of $(Q,f)$ up to $S\-$conjugation implies that
 $$\hat H_\theta = \hat S\.N_{\hat H_\theta} (Q,f)\qq  A_\theta \cong  N_{\hat H_\theta} (Q,f)/N_{\hat S} (Q,f)
\eqno £14.1.$$
On the other hand, since $A_\theta$ is a {\it cyclic\/} $p'\-$group, considering the $k^*\-$group 
$C_{\hat S}^{\hat H_\theta}(Q,f)$ defined in~[3,~15.5.4], it follows from [3,~Lemma~15.16] that $(Q,f)$ 
splits into a set $\{(Q,f_\rho)\}_\rho$ of selfcentralizing Brauer $\hat H_\theta\-$pairs,
where $\rho$ runs over the set of $k^*\-$sections ${\rm Hom}_{k^*}(C_{\hat S}^{\hat H_\theta}(Q,f),k^*)\,,$ 
and then it is quite clear that $\F_{\!(d,\hat S)} (Q,f)$ is a normal subgroup of $\F_{\!(c_\rho,\hat H_\theta)} (Q,f_\rho)$
where $c_\rho$ denotes the block of $\hat H_\theta$ determined by $(Q,f_\rho)\,.$

\medskip
£15. Consequently, the obvious direct summand $\G_k \big(\hat\F_{\!(c_\rho,\hat H_\theta)} (Q,f_\rho)\!\mid\! \bar \theta^*\big)$ of $\G_k \big(\hat\F_{\!(c_\rho,\hat H_\theta)} (Q,f_\rho)\big)$ clearly corresponds to blocks of
$\hat\F_{\!(c_\rho,\hat H_\theta)} (Q,f_\rho)$ of defect zero and, as above, it is not difficult to prove from  {\it Clifford theory\/} that we have a canonical isomorphism 
$$\sum_{\rho} \G_k \big(\hat\F_{\!(c_\rho,\hat H_\theta)} (Q,f_\rho)\!\mid\! \bar \theta^*\big) 
\cong \G_k (\hat A^{^{\bar\theta^*}}_\theta)
\eqno £15.1\phantom{.}$$
where $\rho$ runs over ${\rm Hom}_{k^*}(C_{\hat S}^{\hat H_\theta}(Q,f),k^*)$  
and $\hat A^{^{\bar\theta^*}}_\theta$ is a suitable central $k^*\-$ex-tension  of $A_\theta\,;$ once again, since $A_\theta$
 is cyclic, the  $k^*\-$extension $\hat A^{^{\bar\theta^*}}_\theta$ is split and, since ${\rm Aut} (A_\theta)$ is abelian,
this automorphism group acts {\it canonically\/} on~$\hat A^{^{\bar\theta^*}}_\theta\,.$

\medskip
£16. Finally, for any irreducible Brauer character $\theta$ of $\hat S$ in the block $d\,,$
 from isomorphisms~£12.2 and~£15.1 we get an isomorphism
 $$\G_k (\hat H_\theta\!\mid\! \theta)\cong \sum_{\rho} \G_k \big(\hat\F_{\!(c_\rho,\hat H)} 
 (Q,f_\rho)\!\mid\! \bar \theta^*\big) 
 \eqno £16.1\phantom{.}$$
where $\rho$ runs over ${\rm Hom}_{k^*}(C_{\hat S}^{\hat H_\theta}(Q,f),k^*)\,,$  which is compatible with 
the action of ${\rm Aut} (A_\theta)\,,$ and therefore with the action  of~$C_\theta\,;$ moreover, it is easily checked that
it is compatible with the action of the blocks of $\hat H\,;$ hence, we get
$${\rm rank}_\O \big(\G_k (\hat H_\theta,c_\theta\!\mid\! \theta)^{C_\theta}\big) = 
\sum_{\rho} \G_k \big(\hat\F_{\!(c_\rho,\hat H_\theta)} (Q,f_\rho)\!\mid\! \bar \theta^*\big) 
\eqno £16.2\phantom{.}$$
where $\rho$ runs over the elements of ${\rm Hom}_{k^*}(C_{\hat S}^{\hat H_\theta}(Q,f),k^*)$ such that
$(Q,f_\rho)$ is a Brauer $(c_\theta,\hat H_\theta)\-$pair. Now, according to equality~£12.1, the equality in~£10.1 clearly follows from the corresponding sum of these equalities.

\medskip
£17. In conclusion, in order to get Alperin's Conjecture, it suffices to verify that, {\it for  any block $(c,\hat H)$ 
having a normal sub-block~$(d,\hat S)$ of positive defect  such that the $k^*\-$quotient $S$ of $\hat S$ is simple, 
$H/S$ is a cyclic $p'\-$group and $C_{H}(S) = \{1\}\,,$ and for any cyclic subgroup $C$ of ${\rm Out}_{k^*}(\hat H)_c$ we have
$${\rm rank}_\O\big(\G_k (\hat H,c)^C\big) = \sum_{(Q,f)} {\rm rank}_{\O} \Big(\G_k \big(\hat\F_{\!(c,\hat H)} 
(Q,f),\bar b_f\big)^{C_f}\Big)
\eqno £17.1\phantom{.}$$
where $(Q,f)$ runs over a set of representatives for the set of $H\-$conjugacy classes of selfcentralizing 
Brauer $(c,\hat H)\-$pairs and,  for such a pair, $\bar b_f$ denotes the sum of blocks of defect zero of 
$\hat\F_{\!(c,\hat H)} (Q,f)/\F_Q (Q)\,,$ and  $C_f$ the ``stabilizer'' of~$(Q,f)$ in $C\,.$\/} We honestly believe that
this condition is really easier than the checking demanded in~[2]. Of course, this condition could be true whereas  the statement (Q) in [4,~1.4] failed!

\bigskip
\bigskip
\noindent
{\bf References}
\bigskip
\noindent
[1]\phantom{.} Reinhard Kn\"orr and Geoffrey Robinson, {\it Some remarks on a
conjecture of Alperin}, Journal of London Math. Soc. 39(1989), 48-60.
\smallskip\noindent
[2]\phantom{.} Gabriel Navarro and Pham Huu Tiep, {\it ``A reduction theorem for the Alperin weight conjecture''\/},
to appear in Inventiones math.
\smallskip\noindent
[3]\phantom{.} Llu\'\i s Puig, {\it ``Frobenius categories versus Brauer blocks''\/}, Progress in Math. 
274(2009), Birkh\"auser, Basel.
\smallskip\noindent
[4]\phantom{.} Llu\'\i s Puig, {\it On the reduction of Alperin's Conjecture to the quasi-simple groups\/}, Journal of Algebra, 
328(2011), 372-398

\vfill

\end